\documentclass[runningheads]{tex/CPAIOR/proceedings_template/llncs}
\usepackage{tabulary}
\usepackage{amsmath}
\usepackage{amssymb}
\usepackage{todonotes}
\usepackage[T1]{fontenc}
\usepackage{subfigure}
\usepackage{graphicx}
\usepackage{float}
\allowdisplaybreaks

\usepackage[colorlinks = true,
            linkcolor = blue,
            urlcolor  = blue,
            citecolor = red,
            anchorcolor = blue]{hyperref}

\begin{document}
\title{Leveraging Quantum Computing for \\
Accelerated Classical Algorithms \\
in Power Systems Optimization}
\titlerunning{Leveraging Quantum Computing for Power Systems Optimization}
%
\author{Rosemary Barrass\inst{1}\orcidID{0009-0007-9858-1876} \and
Harsha Nagarajan\inst{2}\orcidID{0000-0003-4550-1100} \and
Carleton Coffrin\inst{2}\orcidID{0000-0003-3238-1699}}
\authorrunning{R. Barrass et al.}
%
\institute{Industrial \& Systems Engineering, Georgia Institute of Technology, GA, USA \and 
Los Alamos National Laboratory, Los Alamos, NM, USA \\
\href{mailto:rbarrass3@gatech.edu}{rbarrass3@gatech.edu}, \href{mailto:harsha@lanl.gov}{harsha@lanl.gov}, \href{mailto:cjc@lanl.gov}{cjc@lanl.gov}}
\maketitle              
\begin{abstract}
The recent advent of commercially available quantum annealing hardware (QAH) has expanded opportunities for research into quantum annealing-based algorithms. 
In the domain of power systems, this advancement has driven increased interest in applying such algorithms to mixed-integer problems (MIP) like Unit Commitment (UC).
UC focuses on minimizing power generator operating costs while adhering to physical system constraints. 
Grid operators solve UC instances daily to meet power demand and ensure safe grid operations.
This work presents a novel hybrid algorithm that leverages quantum and classical computing to solve UC more efficiently. 
We introduce a novel Benders-cut generation technique for UC, thereby enhancing cut quality, reducing expensive quantum-classical hardware interactions, and lowering qubit requirements.
Additionally, we incorporate a $k$-local neighborhood search technique as a recovery step to ensure a higher quality solution than current QAH alone can achieve. 
The proposed algorithm, QC4UC, is evaluated on a modified instance of the IEEE RTS-96 test system.
Results from both a simulated annealer and real QAH are compared, demonstrating the effectiveness of this algorithm in reducing qubit requirements and producing near-optimal solutions on noisy QAH.

\end{abstract}
\section{Introduction}
We stand on the edge of a technological revolution, where the rise of quantum computing promises to redefine our capabilities across a myriad of fields, unlocking new pathways for solving a wide range of complex optimization problems. 
Notably, quantum computers have demonstrated great promise in tackling combinatorial optimization problems, offering a glimpse into the potential for quantum computing to reshape complex decision-making processes.
One area poised for significant transformation is power systems optimization, where mixed-integer optimization problems like unit commitment (UC) must be solved multiple times every day to maintain safe operations of massive power grids. 
UC is a fundamental optimization problem in power grid operations that determines which generating units should be online at any given time to meet anticipated electricity demand over a specified planning horizon, typically ranging from hours to days. 
The goal of UC is to minimize overall energy production costs while satisfying engineering constraints, such as generation limits, minimum up and down times of power plants, and operational reliability requirements. 
Effective UC decision-making is crucial for maintaining the stability and efficiency of power systems, as it ensures that there is a sufficient supply of electricity to meet consumer needs while minimizing wastage of resources. 
For this reason, researchers have started to explore the potential of quantum and hybrid quantum-classical (HQC) algorithms for speeding up the solution of UC and other critical power systems optimization problems.

Research on quantum and HQC algorithms for power systems optimization remains in its early stages, hindered by the limited availability of QPUs and the inherent constraints of NISQ hardware. Demonstrating quantum advantage in practical applications remains challenging due to three fundamental issues. First, QPUs suffer from noise and decoherence, leading to unreliable computations \cite{QCTakingOnItsBiggestChllnge:Noise2024,QCInNISQEraAndBynd2018}. Second, their small scale restricts computational capacity, delaying the realization of quantum advantage. While IBM and Google have recently surpassed the 1000-qubit threshold, the effective number of programmable qubits remains significantly lower due to the overhead of error correction \cite{RcrdBrkingQCHasMoreThn1000Qbts2023}. Third, quantum circuit generation and embedding introduce substantial time complexity, with D-Wave’s own embedding algorithm exhibiting exponential scaling in practical implementations \cite{APrctclHeurstcForFndngGrphMnrs2014}. These challenges underscore the pressing need for algorithmic strategies that can adapt to noisy outputs, limited qubit resources, and high computational overhead.

Despite the limitations of NISQ-era QPUs, evidence suggests that certain power system problems are well-suited for testing on small, noisy quantum hardware \cite{BrdgingGaptoNxtGenPSPlningAndOpswithQC2024}. 
HQC algorithms have shown promise as an intermediary approach for solving MIPs with existing quantum technology. Hybrid quantum-classical algorithms have been developed for a range of mixed-integer power system applications, including unit commitment \cite{OnHQCAlgsForMIP2022,feng2022novel,AdptingQAOAforUC2021,HQCMultCutBDApprWPSApp2023,ASclblFllyDistrbtedQADMMforUC2024,FstQAforSrchQuasOptSolsofUC2024}, micro-grid scheduling \cite{HQCGnrlBDAlgforUCwithMultNtwrkMcrogrds2022,QuanDistrUC:AppInMcroGrds2022}, multi-energy system optimization \cite{IntegQandCCforMultEngySysOptusingBD2024}, and optimal transmission switching \cite{HQCAlgForMIPinPS2024}. Among HQC-based unit commitment approaches, research is largely divided between quantum approximation optimization algorithms (QAOA) and quantum annealing (QA). Koretsky et al. applied QAOA within a hybrid framework, partitioning UC into an inner binary optimization problem solved by a QPU and an outer linear optimization problem handled classically \cite{AdptingQAOAforUC2021}. This formulation, common in HQC algorithms, reduces qubit requirements for large-scale problems, enabling demonstrations of quantum speed-up on current hardware. In Koretsky et al.'s study, this speed-up is observed; however, due to hardware limitations, the quantum circuit is simulated using simulated annealing (SA) rather than executed on a physical QPU.


Paterakis introduced an HQC algorithm for solving the UC problem, leveraging a QPU to optimize cut selection within a Benders' decomposition framework \cite{Benders1962,HQCMultCutBDApprWPSApp2023}. The algorithm demonstrated the benefits of quantum-assisted optimal cut selection on an 8-bus system. However, the constraints of NISQ hardware led to excessive embedding overhead, diminishing its viability compared to purely classical approaches. Notably, Paterakis' approach struggled with scalability, failing to improve performance on a larger 30-bus system--despite the expectation that quantum advantage should manifest as problem size grows. This study underscores the inherent challenges of current noisy QPUs, which remain ill-suited for handling large-scale combinatorial optimization in their present form.

These algorithms highlight the necessity for designing HQC methods that account for NISQ-era constraints, including noise-induced errors in large problems, qubit embedding challenges, and excessive iteration overhead from quantum-classical data transfer. This paper presents QC4UC, a novel HQC algorithm tailored for the UC problem, explicitly addressing the limitations of small-scale noisy quantum hardware. By harnessing quantum computing's unique capabilities alongside classical optimization, QC4UC enhances solution accuracy while reducing qubit requirements compared to previous hybrid approaches.


\section{Methodology}
\begin{figure}
    \hspace{-1.5cm}
    \includegraphics[width=1.2\linewidth]{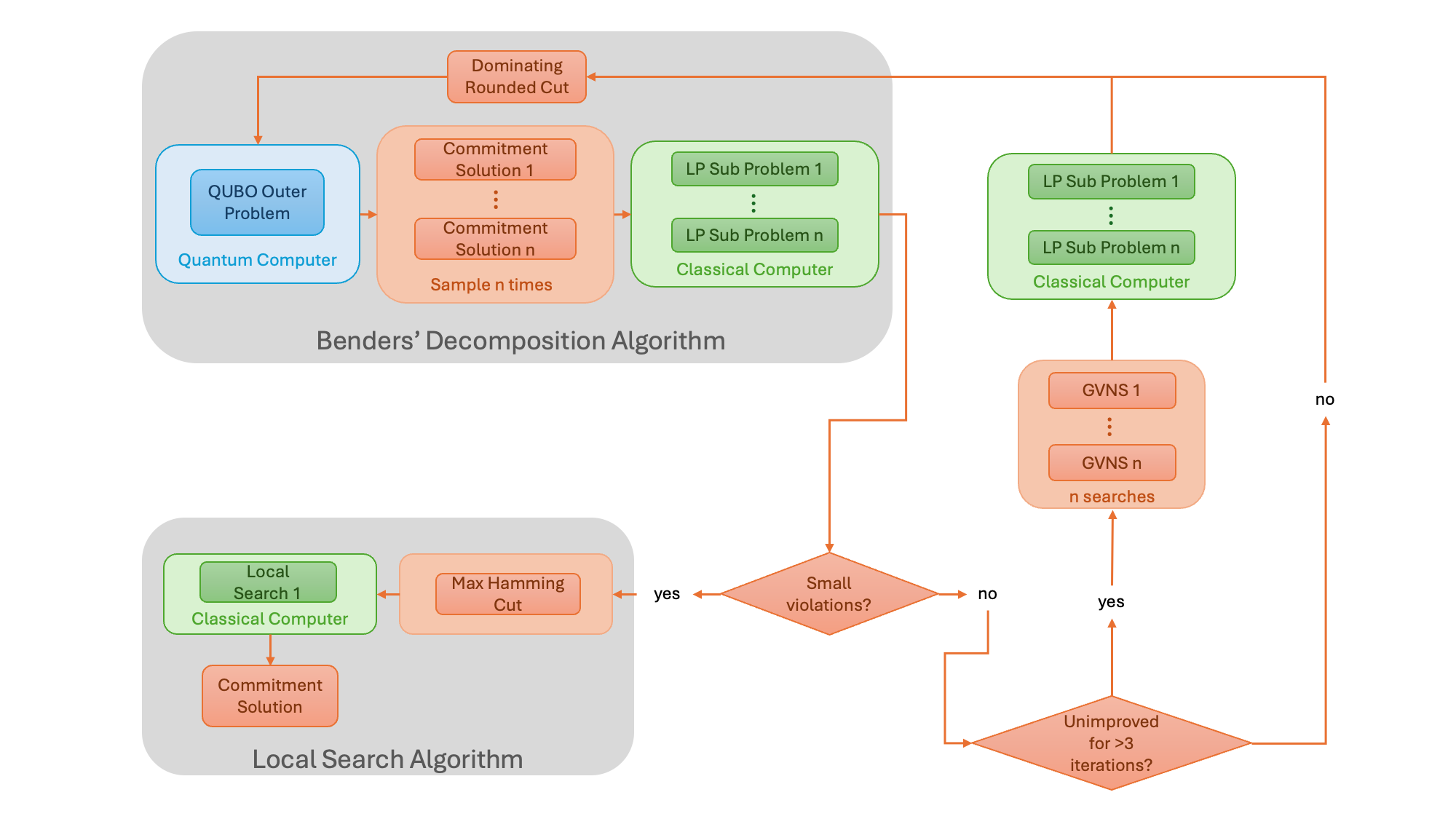}
    \caption{The proposed iterative algorithm with rounded cuts, general variable neighborhood search step, and $k$-local neighborhood search recovery step.}
    \label{sec:methodology:fig:algorithm}
\end{figure}
This paper presents a novel hybrid algorithm for operation on NISQ-era quantum computers which implements a recovery step for smoothing out noise and clever cut generation to reduce both qubit requirements and quantum-classical interactions.
The proposed algorithm is comprised of two sub-algorithms:
\begin{itemize}
    \item A classical Benders' decomposition algorithm with a QUBO main problem (MP) and linear programming (LP) sub-problems (SP).
    The QUBO MP is solved using either SA or QA, gathering $n$ samples to at each iteration to be used in $n$ SP.
    In this paper, the QUBO is solved using both SA and QA for comparison.
    The $n$ SP are solved on a classical computer to find the largest violating QUBO solution.
    The SP with the largest violating solution is then used to generate a rounded cut, which is added to the QUBO.
    This algorithm is repeated until the load and generation penalties are small, with some upper bound set by the user.
    If, during the execution of this algorithm, some $i$ iterations do not show large improvements in the lower bound, a general variable neighborhood search (GVNS) is employed \cite{VarNghbrhdSrch1997}.
    This step effectively ``shakes'' the solution if a local minimum has been found that is difficult to escape with annealing techniques.
    This solution shaking acts as a way to ``reset'' the current solution without requiring step size tuning in annealing algorithm.
    \item A $k$-local neighborhood search \cite{LocalBranching2003} is employed after the completion of the Benders' decomposition algorithm to act as a feasibility recovery step.
    The user picks some $k$ integer value, which should be small.
    This $k$ is the upper bound for a maximum Hamming distance cut, which is added to the original full MIP formulation.
    This cut greatly reduces the search space for the MIP and acts as a recovery step to bring poor QA solutions to the feasible region of the original problem.
\end{itemize}
This new algorithm is applied to a realistic UC formulation, which reflects the pertinent power systems and engineering constraints. 
The formulation includes DC power flow constraints and thermal line limits, the latter of which is not included in prior HQC literature. 
These constraints, however, are standard in power systems and represent critical operating constraints in real-world systems. 
The inclusion of thermal line constraints introduces additional complexity to the MIP, which our proposed QC4UC algorithm can successfully solve for larger instances than previous work. 
We compare our QC4UC algorithm with a classical version of QC4UC which uses SA to solve the QUBO MP.

The full UC formulation is described in Section \ref{sec:formulation}. 
The Benders' decomposition algorithm with GVNS and dynamic binary encoding is described in Section \ref{sec:BDalgorithm}. 
The $k$-local neighborhood search recovery step is described in Section \ref{sec:KLocalSearch}.
Computational results are presented in Section \ref{sec:results}.

\subsection{Problem Formulation}\label{sec:formulation}

The fundamental UC formulation, which serves as the basis for our approach, is given below. In this formulation, $\mathcal{T} = \{1, 2, \dots, T\}$ represents the set of time steps $t$ from $1$ to $T$, and $\mathcal{G} = \{1, 2, \dots, G\}$ denotes the set of generators $g$ from $1$ to $G$. All variables are bolded to differentiate them from constants. Each generator $g$ at time $t$ is characterized by three binary variables: $\mathbf{u}_{g,t}$ (on/off), $\mathbf{v}_{g,t}$ (startup), and $\mathbf{w}_{g,t}$ (shutdown). A generator must remain online for at least $T^{\text{minup}}_g$ time steps after startup and offline for at least $T^{\text{mindn}}_g$ after shutdown. Marginal generator dispatch—representing energy dispatched beyond the minimum required for an online generator—is denoted by the continuous variable $\mathbf{p}_{g,t}$. This value is constrained between $0$ and $\overline{p}_{g,t} - \underline{p}_{g,t}$, the difference between maximum and minimum dispatch levels. At each time step $t$, the sum of all marginal and minimum generator dispatch must satisfy the demand $d_t$. Each generator incurs a minimum running cost $c_1$, a startup cost $c_2$, a shutdown cost $c_3$, and quadratic dispatch costs $c_4$ and $c_5$.
\begin{align}
    \min_{\mathbf{u}, \mathbf{p}} \quad & c_1^T \mathbf{u} + c_2^T\mathbf{v} + c_3^T\mathbf{w} + (c_4)^T (\mathbf{p} \otimes \mathbf{p}) + (c_5)^T \mathbf{p} \quad \quad \tag{\text{Basic UC}}\label{sec:formulation:eq:basic-uc}\\ 
    \text{s.t. } \quad &  \mathbf{u}_{g, t} - \mathbf{u}_{g, t-1} = \mathbf{v}_{g, t} - \mathbf{w}_{g, t}, \quad \quad \quad \forall g,t \in \mathcal{G}, \mathcal{T} \label{sec:formulation:eq:basic-uc:susd-logic}\\
    \quad & \sum_{\tau = t - T^{\text{minup}}_g + 1}^{t} \mathbf{v}_{g, \tau} \leq \mathbf{u}_{g, t}, \quad \quad \quad \forall g,t \in \mathcal{G}, \{T^{\text{minup}}_g, \dots, \mathcal{T}\} \label{sec:formulation:eq:basic-uc:minup}\\
    \quad & \sum_{\tau = t - T^{\text{mindn}}_g + 1}^{t} \mathbf{w}_{g, \tau} \leq 1 - \mathbf{u}_{g, t}, \quad \quad \forall g,t \in \mathcal{G}, \{T^{\text{mindn}}_g, \dots, \mathcal{T}\} \label{sec:formulation:eq:basic-uc:mindn}\\
    \quad & \sum_{g \in \mathcal{G}} \mathbf{p}_{g,t} + \sum_{g \in \mathcal{G}} \left( \underline{p}_{g,t}\cdot\mathbf{u}_{g, t} \right) = d_t, \quad \quad \forall t \in \mathcal{T} \label{sec:formulation:eq:basic-uc:pb}\\
    \quad & \mathbf{u}_{g, t}, \ \mathbf{v}_{g, t}, \ \mathbf{w}_{g, t} \in \{0, 1\}, \quad \quad \quad \forall g, t \in \mathcal{G}, \mathcal{T} \label{sec:formulation:eq:basic-uc:u-bin}\\
    \quad & 0 \leq \mathbf{p}_{g, t} \leq  \overline{p}_{g,t} - \underline{p}_{g,t}, \quad \quad \quad \forall g, t \in \mathcal{G}, \mathcal{T} \label{sec:formulation:eq:basic-uc:pmin-pmax}
\end{align}
In this paper, costs are linearized using a piecewise approximation. The generator dispatch variables for each generator $g$ are partitioned into $L$ segments between the minimum and maximum dispatch levels.  The fraction of the total energy dispatch allocated to segment $l$ for generator $g$ at time $t$ is denoted by $\boldsymbol{\alpha}^l_{g,t}$. The cost coefficient for segment $l$ in the piecewise linear approximation of the generation cost is given by ${c_6}^l_{g,t}$. Thus, the total marginal energy dispatch for generator $g$ at time $t$ is the sum of these segments:
\begin{align}
    \label{sec:formulation:eq:marginal-dispatch-eq}
    \mathbf{p}_{g, t} = \sum_{l \in \mathcal{L}_g}({p}^l_{g, t} - \underline{p}_{g, t})\cdot\boldsymbol{\alpha}^l_{g,t}, & \quad \forall g, t \in \mathcal{G}, \mathcal{T} 
\end{align}
All of these fractions must add up to $1$ when the generator is online and $0$ when the generator is offline:
\begin{equation}
    \begin{aligned}
        \label{sec:formulation:eq:fractional-eq}
        \sum_{l \in \mathcal{L}_g} \boldsymbol{\alpha}^l_{g, t} = \mathbf{u}_{g, t}, & \quad \forall g, t \in \mathcal{G}, \mathcal{T} 
    \end{aligned}
\end{equation}
In addition to these linearizing constraints, practical engineering constraints are incorporated to enhance solution quality. We define $ B_{n,j} $ as the susceptance of the power line connecting bus $ n $ to bus $ j $, and $ \boldsymbol{\theta}_{n,t} $ as the phase angle at bus $ n $ at time $ t $. The set $ \mathcal{N} = \{1, \dots, N\} $ represents all buses in the network. To ensure that the thermal limits $ F^{\text{max}}_{n,j} $ of each line $(n,j)$ are not exceeded, we impose the following constraint:
\begin{equation}
    \begin{aligned}
        \label{sec:formulation:eq:line-constraint}
        B_{n, j}|(\boldsymbol{\theta}_{n, t} - \boldsymbol{\theta}_{j, t})| \leq F^{\text{max}}_{n, j}, & \quad  \forall n \in \mathcal{N} 
    \end{aligned}
\end{equation}
The following ramping constraints are introduced to limit the change in power dispatch that each generator can achieve between consecutive time steps. $ R^{\text{startup}}_g $ and $ R^{\text{up}}_g $ define the maximum allowable increases in energy dispatch during startup and between successive time steps, respectively. Similarly, $ R^{\text{shutdown}}_g $ and $ R^{\text{down}}_g $ specify the maximum allowable decreases during shutdown and between consecutive time steps. The variable $ \mathbf{r}_{g,t} $ represents the spinning reserve provided by generator $ g $ at time $ t $.
\begin{align}
    \mathbf{p}_{g, t} + \mathbf{r}_{g, t}  \leq (\overline{p}_{g, t} - \underline{p}_{g, t})\cdot\mathbf{u}_{g, t} \quad &\quad \nonumber \\
     \quad - \max\{\overline{p}_{g, t} - R^{\text{startup}}_g, 0\}\cdot\mathbf{v}_{g, t}, & \quad  \forall g, t \in \mathcal{G}, \mathcal{T} \label{sec:formulation:eq:ramping-constraints:suramp} \\
    \mathbf{p}_{g, t} + \mathbf{r}_{g, t}  \leq (\overline{p}_{g, t} - \underline{p}_{g, t})\cdot\mathbf{u}_{g, t} \quad &\quad \nonumber \\
     \quad - \max\{\overline{p}_{g, t} - R^{\text{shutdown}}_g, 0\}\cdot\mathbf{w}_{g, t+1}, & \quad \forall g, t \in \mathcal{G}, \mathcal{T}\backslash \{T\} \label{sec:formulation:eq:ramping-constraints:sdramp}\\
    \mathbf{p}_{g, t} + \mathbf{r}_{g, t} - \mathbf{p}_{g, t-1} \leq R^{\text{up}}_g, & \quad \forall g, t \in \mathcal{G}, \mathcal{T} \backslash \{1\} \label{sec:formulation:eq:ramping-constraints:rampup}\\
    \mathbf{p}_{g, t-1} - \mathbf{p}_{g, t} \leq R^{\text{down}}_g, & \quad \forall g, t \in \mathcal{G}, \mathcal{T} \backslash \{1\} \label{sec:formulation:eq:ramping-constraints:rampdn}
\end{align}
Including constraints $\eqref{sec:formulation:eq:marginal-dispatch-eq}$
through \eqref{sec:formulation:eq:ramping-constraints:rampdn} to $(\ref{sec:formulation:eq:basic-uc})$ gives us the final MILP UC formulation.
$\underline{c_6}_{g,t}$ is added to represent the minimum generation dispatch cost for generator $g$ at time $t$ (${c_6}^1_{g, t} = \underline{c_6}_{g,t}$).
\begin{align}
    \min_{\mathbf{u}, \mathbf{p}} \quad & c_1^T \mathbf{u} + c_2^T\mathbf{v} + c_3^T\mathbf{w} + \sum_{l \in \mathcal{L}_g} ({c_6}^l_{g, t} - \underline{c_6}_{g,t}) \boldsymbol{\alpha}^l_{g,t} & \tag{\text{Full UC}}\label{sec:formulation:eq:full-uc} \\
    \text{s.t. } \quad &  \eqref{sec:formulation:eq:basic-uc:susd-logic}, \eqref{sec:formulation:eq:basic-uc:minup}, \eqref{sec:formulation:eq:basic-uc:mindn}, \eqref{sec:formulation:eq:basic-uc:pmin-pmax}, \eqref{sec:formulation:eq:marginal-dispatch-eq}, \eqref{sec:formulation:eq:fractional-eq}, \eqref{sec:formulation:eq:line-constraint}, \eqref{sec:formulation:eq:ramping-constraints:suramp}, \eqref{sec:formulation:eq:ramping-constraints:sdramp}, \eqref{sec:formulation:eq:ramping-constraints:rampup}, \eqref{sec:formulation:eq:ramping-constraints:rampdn} \nonumber\\
    \quad & \sum_{g \in \mathcal{G}} \mathbf{p}_{g,t} + \sum_{g \in \mathcal{G}} \left( \underline{p}_{g,t}\cdot\mathbf{u}_{g, t} \right) - d_{n, t} = &\nonumber \\
    \quad & \quad \sum_{(n, j) \in \mathcal{L}} B_{n, j} (\boldsymbol{\theta}_{n, t} - \boldsymbol{\theta}_{j, t}),& \forall n, t \in \mathcal{N}, \mathcal{T} \label{sec:formulation:eq:full-uc:pb}\\
    \quad & \mathbf{u}_{g, t}, \ \mathbf{v}_{g, t}, \ \mathbf{w}_{g, t}  \in \{0, 1\}, & \forall g, t \in \mathcal{G}, \mathcal{T} \label{sec:formulation:eq:full-uc:bin-vars}\\
    \quad & \mathbf{r}_{g, t} \geq 0, & \forall g, t \in \mathcal{G}, \mathcal{T} \label{sec:formulation:eq:full-uc:ramp-var}\\
    \quad & \boldsymbol{\alpha}^l_{g,t} \geq 0, & \forall l, g, t \in \mathcal{L}_g, \mathcal{G}, \mathcal{T} \label{sec:formulation:eq:full-uc:piecewise-disp-var}\\
    \quad & \boldsymbol{\theta}_{\text{ref}, t} = 0, & \forall t \in \mathcal{T} \label{sec:formulation:eq:full-uc:volt-ang-var}
\end{align}
\subsection{Benders' Decomposition Algorithm}\label{sec:BDalgorithm}
This section presents the proposed Benders' decomposition algorithm (QC4UC), incorporating a novel combination of a dynamic precision encoding technique, rounded cut generation, and a general variable neighborhood search. The MILP formulation \eqref{sec:formulation:eq:full-uc} is decomposed into a binary optimization MP, solved using QAH, and a continuous SP, handled by a commercial branch-and-bound solver on a CPU. $\boldsymbol{\eta}$ represents the MP's estimation of the SP objective value. The MP formulation is as follows:
\begin{align}
    \min_{\mathbf{u}, \mathbf{p}} \quad & c_1^T \mathbf{u} + c_2^T\mathbf{v} + c_3^T\mathbf{w} + \boldsymbol{\eta} & \tag{\text{MP}}\label{sec:BDalgorithm:eq:MP}\\
    \text{s.t. } \quad &  \eqref{sec:formulation:eq:basic-uc:susd-logic}, \eqref{sec:formulation:eq:basic-uc:minup}, \eqref{sec:formulation:eq:basic-uc:mindn}, \eqref{sec:formulation:eq:full-uc:bin-vars} \nonumber
\end{align}
The MP generates solutions ${u}^*, {v}^*, {w}^*$ which represent the commitment, startup, and shutdown decisions.
In order to ensure all MP solutions are feasible in the subproblem, $\boldsymbol{\delta}^+_t, \boldsymbol{\delta}^-_t$ are added to the SP as penalty terms for relaxing the power balance, where $\boldsymbol{\delta}^+_t$ represents overproduction of energy and $\boldsymbol{\delta}^-_t$ represents underproduction of energy.
The SP take the form:
\begin{align}
    Q({u}^*) = \quad & \min_{\boldsymbol{\alpha}, \mathbf{p}} \sum_{l \in \mathcal{L}_g} ({c_6}^l_{g, t} - \underline{c_6}_{g,t}) \boldsymbol{\alpha}^l_{g,t} + c^{\text{pen}}_t (\boldsymbol{\delta}^+_t - \boldsymbol{\delta}^-_t) & \tag{\text{SP}}\label{sec:BDalgorithm:eq:SP}\\
    \text{s.t. } & \eqref{sec:formulation:eq:basic-uc:pmin-pmax}, \eqref{sec:formulation:eq:marginal-dispatch-eq}, \eqref{sec:formulation:eq:fractional-eq}, \eqref{sec:formulation:eq:line-constraint}, \eqref{sec:formulation:eq:ramping-constraints:rampup}, \eqref{sec:formulation:eq:ramping-constraints:rampdn}, \eqref{sec:formulation:eq:full-uc:ramp-var}, \eqref{sec:formulation:eq:full-uc:piecewise-disp-var}, \eqref{sec:formulation:eq:full-uc:volt-ang-var} & \nonumber \\
    \quad & \sum_{g \in \mathcal{G}} \mathbf{p}_{g,t} + \sum_{g \in \mathcal{G}} \left( \underline{p}_{g,t}\cdot{u}^*_{g, t} \right) - d_{n, t} = &\nonumber \\
    \quad & \quad \sum_{(n, j) \in \mathcal{L}} B_{n, j} (\boldsymbol{\theta}_{n, t} - \boldsymbol{\theta}_{j, t}) + (\boldsymbol{\delta}^+_t + \boldsymbol{\delta}^-_t), & \forall n, t \in \mathcal{N}, \mathcal{T} \\
    \quad & \eqref{sec:formulation:eq:ramping-constraints:suramp} \text{ with fixed values $u^*_{g, t}, v^*_{g, t}$} \\
    \quad & \eqref{sec:formulation:eq:ramping-constraints:sdramp} \text{ with fixed values $u^*_{g, t}, w^*_{g, t}$} \\
    \quad & \boldsymbol{\delta}^+_t \geq 0, \ \boldsymbol{\delta}^-_t \leq 0, & \forall t \in \mathcal{T}
\end{align}
These penalty terms have an associated penalty cost $c^{\text{pen}}_t$ which is larger than the highest generator dispatch cost in the system.
$\lambda^*_i$'s are introduced as the dual variables for the constraints of the SP.
Note that $\lambda^*_1$ is always zero, as the first SP constraint has a right-hand side (RHS) value of $0$.
The optimality cut generated by each SP takes the form
\begin{align*}
    (\mathbf{u}_{g, t})^T \lambda_2^* \tag{OptCut}\label{sec:BDalgorithm:eq:opt-cut}\\
    + \left(d_{n, t} - \sum_{g \in \mathcal{G}_n} \left( \underline{p}_{g,t}\cdot\mathbf{u}_{g, t} \right)\right)^T \lambda_3^* \\
    + \left((\overline{p}_{g, t} - \underline{p}_{g, t})\cdot\mathbf{u}_{g, t} - \max\{\overline{p}_{g, t} - R^{\text{startup}}_g, 0\}\cdot\mathbf{v}_{g, t} \right)^T \lambda_4^*\\
    + \left((\overline{p}_{g, t} - \underline{p}_{g, t})\cdot\mathbf{u}_{g, t} - \max\{\overline{p}_{g, t} - R^{\text{shutdown}}_g, 0\}\cdot\mathbf{w}_{g, t+1} \right)^T \lambda_5^* \\
    + (R^{\text{up}}_g)^T \lambda_6^* + (R^{\text{down}}_g)^T \lambda_7^* + (F^{\text{max}}_{n, j})^T \lambda_8^* \\
    \quad \leq \boldsymbol{\eta} 
\end{align*}
In order to solve this problem with QAH, the MP and optimality cuts are converted to a QUBO form.
To do this, we introduce $P_j$, $j = 1, \ldots, 4$ as penalty values that are tuned empirically to penalize distance from the feasible region.
$\mathbf{s}_1, \mathbf{s}_2, \mathbf{s}_3$ are added as slack variables to relax constraints in the MP.
The expression $\text{OptCut}_{\text{LHS}}(\lambda_2^*, \lambda_2^*, \dots, \lambda_8^*)$ is introduced to represent left-hand side (LHS) of the optimality cut with fixed dual variable solutions $\lambda_i^*, i = 2, \dots, 8$ and unfixed binary variables $\mathbf{u}, \mathbf{v},$ and $\mathbf{w}$.
The QUBO form of the MP thus becomes
\begin{align}
    \min_{\mathbf{u}, \mathbf{v}, \mathbf{w}} \quad & c_1^T \mathbf{u} + c_2^T\mathbf{v} + c_3^T\mathbf{w} + \boldsymbol{\eta} & \tag{MP-QUBO}\label{sec:BDalgorithm:eq:MP-QUBO}\\
    \quad & + P_1 \sum_{g \in \mathcal{G}}\sum_{t \in \mathcal{T}} \left(\mathbf{u}_{g, t} - \mathbf{u}_{g, t-1} - \mathbf{v}_{g, t} + \mathbf{w}_{g, t}\right)^2\\
    \quad & + P_2 \sum_{g \in \mathcal{G}}\sum_{t \in \{T^{\text{minup}}_g, \dots, \mathcal{T}\}}\left(\sum_{\tau = t - T^{\text{minup}}_g + 1}^{t} \mathbf{v}_{g, \tau} + \mathbf{s}_1  - \mathbf{u}_{g, t}\right)^2 \\
    \quad & + P_3 \sum_{g \in \mathcal{G}}\sum_{t \in \{T^{\text{mindn}}_g, \dots, \mathcal{T}\}}\left(\sum_{\tau = t - T^{\text{mindn}}_g + 1}^{t} \mathbf{w}_{g, \tau} + \mathbf{s}_2  - 1 + \mathbf{u}_{g, t}\right)^2 \\
    \quad & + P_4 \left( \text{OptCut}_{\text{LHS}}(\lambda_2^*, \lambda_2^*, \dots, \lambda_8^*) - \boldsymbol{\eta} + \mathbf{s}_3 \right)^2 
\end{align}
Note that the slack variables $\mathbf{s}_1, \mathbf{s}_2$ are binary by nature of the inequalities that they are used in.
We can show this by analyzing the minimum up time constraint:
\[\sum_{\tau = t - T^{\text{minup}}_g + 1}^{t} \mathbf{v}_{g, \tau} \leq \mathbf{u}_{g, t}\]
Since that $\mathbf{v}_{g, t}$ and $\mathbf{u}_{g, t}$ are binary for all $g,t \in \mathcal{G}, \mathcal{T}$, $\sum_{\tau = t - T^{\text{minup}}_g + 1}^{t} \mathbf{v}_{g, \tau} - \mathbf{u}_{g, t}$ is bounded between -1 and 0. With slack variable $\mathbf{s}_1$ to create an equality, we have,
\[\sum_{\tau = t - T^{\text{minup}}_g + 1}^{t} \mathbf{v}_{g, \tau} - \mathbf{u}_{g, t} + \mathbf{s}_1 = 0,\]
implying that $\mathbf{s}_1$ can only be a discrete value between 0 and 1. 
A similar argument holds for the minimum down time constraint to confirm its binary nature. 

In contrast, \( \mathbf{s}_3 \) must be discretized for the QUBO formulation.
To address the constraints of NISQ-era quantum hardware, which offers a very limited number of qubits, we implement a novel dynamic precision encoding step and rounded cuts, essential for minimizing the qubits required to binary encode \( \mathbf{s}_3 \) and $\boldsymbol{\eta}$ in each iteration of Benders' decomposition.

\subsubsection{Dynamic Precision Encoding}
The slack variable $\mathbf{s}_3$ can be bounded, which allows for a novel dynamic qubit allocation regime for the cut slack variable at each iteration.
To do this, consider the original optimality cut $\eqref{sec:BDalgorithm:eq:opt-cut}$.
In the MP, the decision variables in this cut will be $\mathbf{u}, \mathbf{v}, \mathbf{w}$ and $\boldsymbol{\eta}$.
For a minimization problem, we can assume that $\boldsymbol{\eta}$ will be chosen such that the left hand side will be equal to $\boldsymbol{\eta}$.
Thus, a worst-case lower bound on the left-hand side can be used to determine an approximate precision necessary for $\boldsymbol{\eta}$.
So, at each iteration of our algorithm, we solve a small optimization problem of the form:
\begin{align*}
    \max_{\mathbf{u}, \mathbf{v}, \mathbf{w}} & \quad \text{OptCut}_{\text{LHS}}(\lambda_2^*, \lambda_2^*, \dots, \lambda_8^*) &\tag{MaxEta}\label{sec:BDalgorithm:eq:dynamic-encoding}\\
    \quad &  \quad \mathbf{u}_{g, t} - \mathbf{u}_{g, t-1} = \mathbf{v}_{g, t} - \mathbf{w}_{g, t}, & \forall g,t \in \mathcal{G}, \mathcal{T} \\
    \quad & \quad \sum_{\tau = t - T^{\text{minup}}_g + 1}^{t} \mathbf{v}_{g, \tau} \leq \mathbf{u}_{g, t}, & \forall g,t \in \mathcal{G}, \{T^{\text{minup}}_g, \dots, \mathcal{T}\} \\
    \quad & \quad \sum_{\tau = t - T^{\text{mindn}}_g + 1}^{t} \mathbf{w}_{g, \tau} \leq 1 - \mathbf{u}_{g, t}, & \forall g,t \in \mathcal{G}, \{T^{\text{mindn}}_g, \dots, \mathcal{T}\}  \\
    \quad & \quad \mathbf{u}_{g, t}, \mathbf{v}_{g, t}, \mathbf{w}_{g, t} \in \{0, 1\}, & \forall g, t \in \mathcal{G}, \mathcal{T} \\
\end{align*}
which is bounded. 
The solution of this IP is then used to inform the precision parameter $p$ to be used in the slack variable encoding for each cut.

\subsubsection{Rounded Optimality Cuts}

In addition to dynamic precision encoding, we introduce a novel rounding technique for the L-shaped cuts of Benders' decomposition to further reduce qubit requirements while ensuring validity as a lower bound to the subproblem’s convex recourse function. Consider the optimality cuts \eqref{sec:BDalgorithm:eq:opt-cut}, simplified to the form:
\[
 a_1^T \mathbf{x} +  a_2^T \mathbf{y} + a_3^T \mathbf{z} + b \leq \boldsymbol{\eta},
\]
where $\boldsymbol{\eta}$ is a decision variable that is minimized in the MP’s objective function to provide a valid lower approximation to the true recourse function of the SP. Here, $ a_1, a_2, a_3 \in \mathbb{R}^n $, $ b \in \mathbb{R} $, and $ \mathbf{x}, \mathbf{y}, \mathbf{z} \in \{0,1\}^n $. To embed this inequality into the QUBO model, we introduce a slack variable $ \mathbf{s}_3 \geq 0 $: $
a_1^T \mathbf{x} +  a_2^T \mathbf{y} + a_3^T \mathbf{z} + b + \mathbf{s}_3 = \boldsymbol{\eta}.
$ Since $ a_1, a_2, a_3, b $ are real numbers, the slack variable $ \mathbf{s}_3 $ is also continuous. To encode $ \mathbf{s}_3 $ with $ p $ decimal places of precision, at least
\[
N = \left\lceil \log_2 \left( \boldsymbol{\eta}^*(10^p) + 1 \right) \right\rceil
\]
binary variables are required, where \( \mathbf{s}_3 \in [0, \eta^*] \) and $\eta^*$ is the objective value of \eqref{sec:BDalgorithm:eq:dynamic-encoding}. 
To reduce the binary variable count, we apply a \textit{conservative rounding} to the coefficients. 
This ensures that the resulting integer-based cut is never more restrictive than the original fractional Benders' cut, thereby preserving a valid lower bound to the convex recourse as \(\boldsymbol{\eta}\) is minimized. Specifically, we define the rounded coefficients \(\tilde{a}_1, \tilde{a}_2, \tilde{a}_3 \in \mathbb{Z}^n\) and \(\tilde{b} \in \mathbb{Z}\), yielding the new constraint:
\[
\tilde{a}_1^T \mathbf{x} + \tilde{a}_2^T \mathbf{y} + \tilde{a}_3^T \mathbf{z} + \tilde{b} + \mathbf{s}_3 = \boldsymbol{\eta}.
\]
To maintain validity, we round each coefficient as follows:   
\begin{itemize}
    \item  Positive coefficients $ ([a_i]_j > 0) $: round down, $ [\tilde{a}_i]_j = \lfloor [a_i]_j \rfloor $. 
    \item Negative coefficients $ ([a_i]_j < 0) $: round up, $ [\tilde{a}_i]_j = \lceil [a_i]_j \rceil $. 
    \item Constant term $(b > 0)$: round down, $\tilde{b} = \lfloor b \rfloor $.
\end{itemize}

With rounded integer coefficients, the slack variable $ \mathbf{s}_3 $ can also be discretized to an integer domain, taking $\eta^* + 1$ possible integer values, requiring only:
\[
M = \left\lceil \log_2 (\eta^* + 1) \right\rceil
\]
binary variables. For large $\eta^*$, the resulting difference in the number of binary variables before and after rounding, $N - M$, simplifies approximately to:
\[
\left\lceil p \log_2 10 \right\rceil \approx \left\lceil 3.322 p \right\rceil.
\]
Thus, rounding the coefficients to the nearest integers and using an integer-valued slack variable significantly reduces the binary variables needed to represent \( \mathbf{s}_3 \), thereby simplifying the QUBO formulation and preserving a valid lower bound, while balancing the trade-off between precision and solution quality.

\subsection{General Variable Neighborhood Search}
QA and SA do not guarantee convergence, meaning the Benders' algorithm may iterate indefinitely without reaching the optimal solution, especially given the reduced precision of the rounded cuts.
This often occurs when the annealer gets trapped in a local minimum that the current step size does not allow ``escape'' from.
One way to address this is through heuristic methods, such as variable flipping or tabu search \cite{TabuSrch1999}, to force solutions out of difficult-to-escape local minima.  
In classical computing, these methods iterate from a single solution, making random modifications to generate new ones. Since they do not guarantee convergence to the global optimum, the final solution from a General Variable Neighborhood Search (GVNS) can vary.
This provides an opportunity to leverage multiple solution samples from the quantum or simulated annealer for improved results. 
The proposed algorithm iterates Benders' cuts until no further solution improvements are observed.
Here, we define some small distance between objective values at each iteration.
If the algorithm produces a solution with an objective value within this distance from the last iteration, then this iteration is considered to be non-improving.
After $i$ consecutive non-improving iterations, it is assumed that the annealer is stuck and the GVNS is employed.
Parallel GVNS will then be performed using $n$ solutions sampled from the quantum or simulated annealer, with $n$ defined by the user. 
Each solution is used to solve the SP, generating a new cut to add to the MP for further iterations of the Benders' decomposition algorithm.

The GVNS includes two stages: a local improvement stage and a shaking stage.
In the local improvement stage, $k_1$ commitment variables are randomly flipped from 0 to 1 or vice versa until a solution that improves the objective value is found. In the shaking stage, $k_2$ commitment variables are randomly flipped without considering the objective value. These stages are applied iteratively across the $n$ QUBO samples.


\subsection{$k$-local Neighborhood Search Algorithm}\label{sec:KLocalSearch}
For the $k$-local neighborhood search algorithm, the full MILP is limited by a maximum Hamming distance constraint which significantly reduces the size of the search space to a ball around the final Benders' decomposition solution.
The maximum Hamming distance constraint takes the form
\begin{align}\label{sec:methodology:eq:minHammingconstraint}
    \quad & \sum_{g, t :  \mathbf{u}^*_{g, t} = 0} \mathbf{u}_{g, t} + \sum_{g, t : \mathbf{u}^*_{g, t} = 1} 1-\mathbf{u}_{g,t} = k + \mathbf{s}_4 
\end{align}
where $\mathbf{s}_4$ is some integer valued slack variable.
This constraint is added to \eqref{sec:formulation:eq:full-uc} and solved on a classical computer.
\section{Computational Results}\label{sec:results}
The proposed algorithm, tested on a significantly larger test case than any previous work, demonstrates effective qubit reduction and achieves near-optimal solutions on noisy QAH.
This is evidenced by comparative results from both SA and real QAH.
In this section, QC4UC is tested on a modified IEEE RTS-96 bus system, which contains 8 thermal generators, 27 renewable generators, and 24 buses.
For the MP, we compare the D-Wave SA package "dwave-neal" version 0.6.0 \cite{DwaveNeal2017} on a classical computer and the D-Wave QA package "dwave-system" version 1.26.0 \cite{DwaveSys2017} on QAH.
The SA package was developed to emulate a noiseless version of D-Wave's current QAH and acts as a way to compare the current NISQ-era QAH to what possible perfect hardware could produce.
The classical computer used to solve the SPs and the SA tests has an Apple M3 Max 16-core CPU.
The QA computations are performed on the Advantage\_system4.1 QPU, which comprises 5612 operational qubits arranged in the Pegasus architecture with a connectivity range of $[-18, 15]$ per qubit, operating at a cryogenic temperature of 15.4 mK. 
The qubit thermal width is 0.198 Ising units, and the quantum critical point is observed at 1.391 GHz, ensuring robust performance during annealing. 
The SPs and recovery step are solved on the classical computer using Gurobi 11.0.2 with default settings. 
A single solution of \eqref{sec:formulation:eq:full-uc} is produced using Gurobi 11.0.2 with default settings on a classical computer to evaluate the solution quality of QC4UC.
It is important to note here that QC4UC is not compared to the Gurobi solution to demonstrate any kind of improvement over classical methods.
Instead, this Gurobi solution is used to demonstrate that QC4UC solutions are good.

\subsection{Solution Quality on Quantum Annealing Hardware}
\begin{figure}[ht]
    \centering
    \subfigure[]{
  \includegraphics[width=0.47\linewidth]{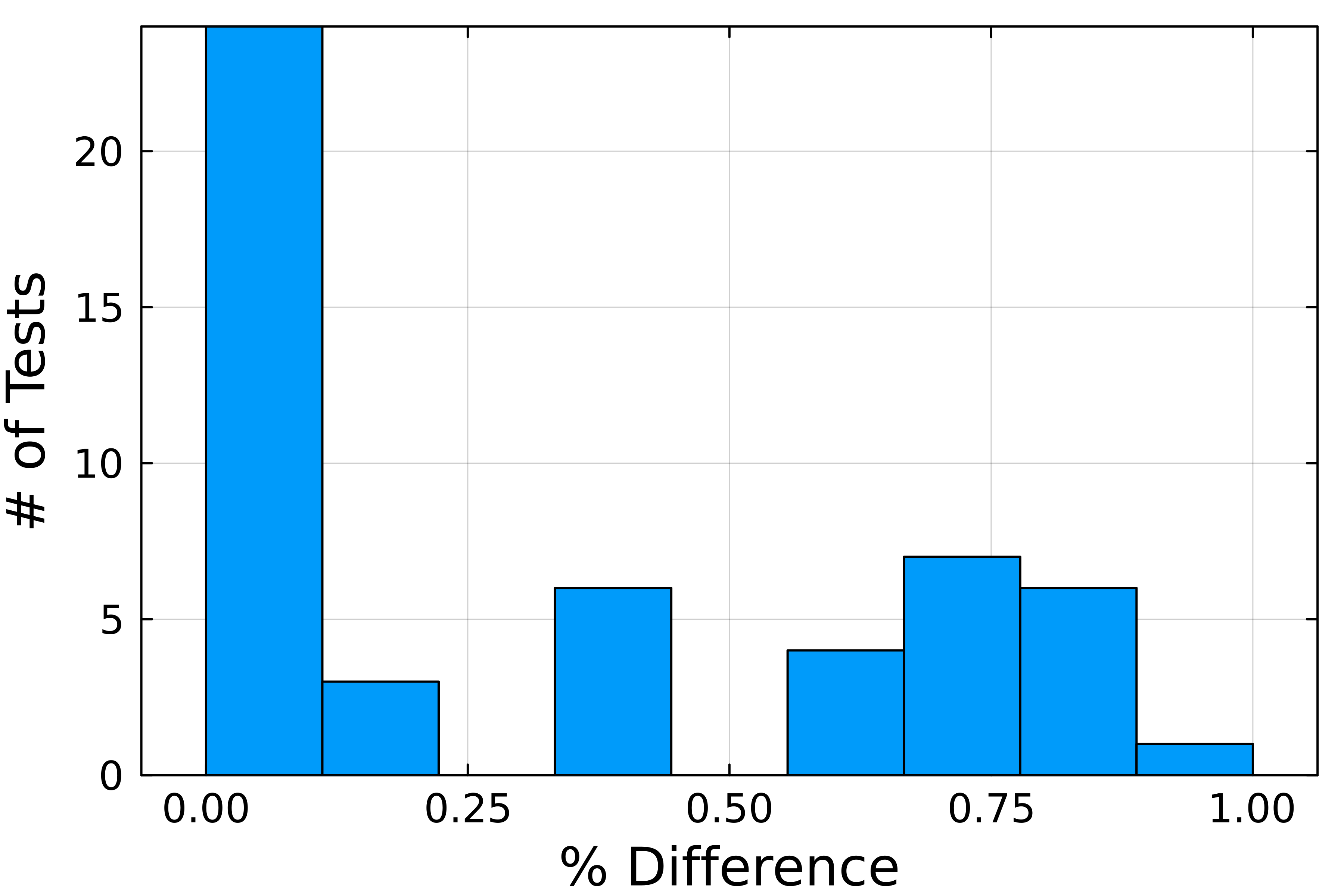}
  }
  \hspace{0.05cm}
    \subfigure[]{
  \includegraphics[width=0.47\linewidth]{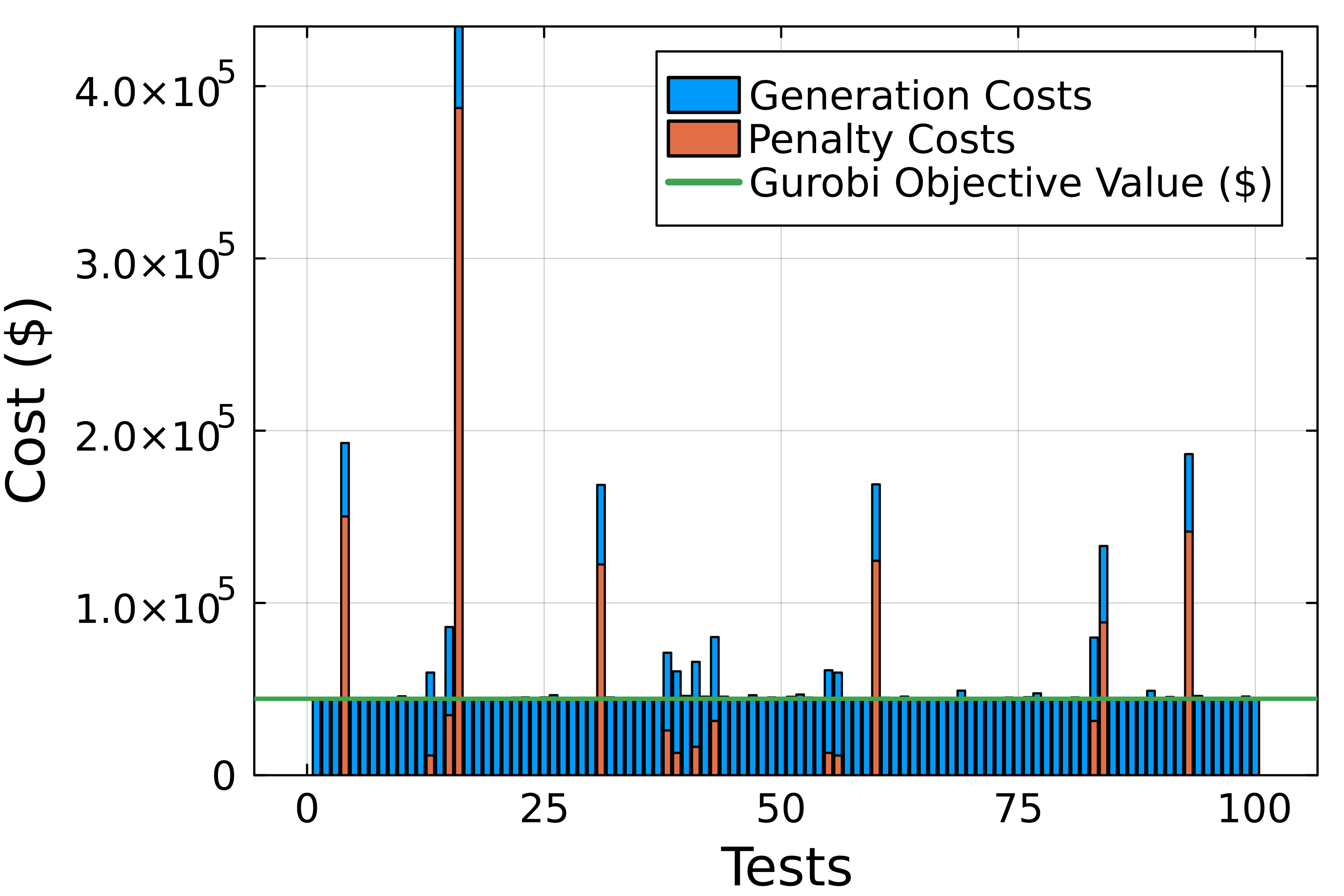}
  }
    \caption{(a) Histogram of the 51 out of 100 QC4UC with QA tests achieving within 1\% of the Gurobi's optimal objective value. (b) Bar plot of penalty and generation costs of all 100 tests of QC4UC with QA.}
    \label{sec:results:fig:percent-obj-val-genpencosts}
\end{figure}

The proposed QC4UC algorithm with QA is solved on the same test case 100 times to assess its consistency and mitigate the impact of noise.
Among the 100 runs, the algorithm produced a solution within 1\% of the Gurobi objective value 51 times. Of these, 24 solutions were within 0.01\% of the Gurobi optimal value, as shown in Figure \ref{sec:results:fig:percent-obj-val-genpencosts}.
This is despite the fact QAH is only guaranteed to converge in probability to a solution by the law of large numbers.
These results demonstrate the efficacy of the algorithm on real NISQ-era quantum hardware.
They also highlight the recovery step's capability to produce high-quality solutions despite significant noise and the potential challenges posed by rounded cuts, such as sub-optimality or cutting into the feasible region of the UC problem. 
Figure \ref{sec:results:fig:percent-obj-val-genpencosts} demonstrates that some of the very large optimality gaps are caused by penalty costs from the load shed penalties.
The penalty costs for load shed are set to $5e5$ \$/MW to act as a big M for keeping load shed to $0$, when possible.
Of the 85 tests with no penalties, 81 of the tests are within 5\% of the Gurobi objective value.
This level of solution quality has yet to be achieved on a problem of this size in the literature.
Overall, this reflects the strength of QC4UC in deriving near-optimal solutions for the UC problem using noisy, stochastic hardware.

Across the same number of independent runs, the QC4UC algorithm with simulated annealing (SA) achieved the optimal objective value for all runs, despite SA's lack of guaranteed convergence properties. 
The same rounded cuts and decomposition methods were applied for these tests as those conducted for QA.
This demonstrates that the proposed algorithm will perform well on noiseless annealing methods.
The consistency of results with SA also highlights the potential for future noiseless QAH.

\subsection{Embedding Algorithm Solve Time and Algorithm Run Time}
\begin{figure}[]
    \centering
    \includegraphics[width=0.8\linewidth]{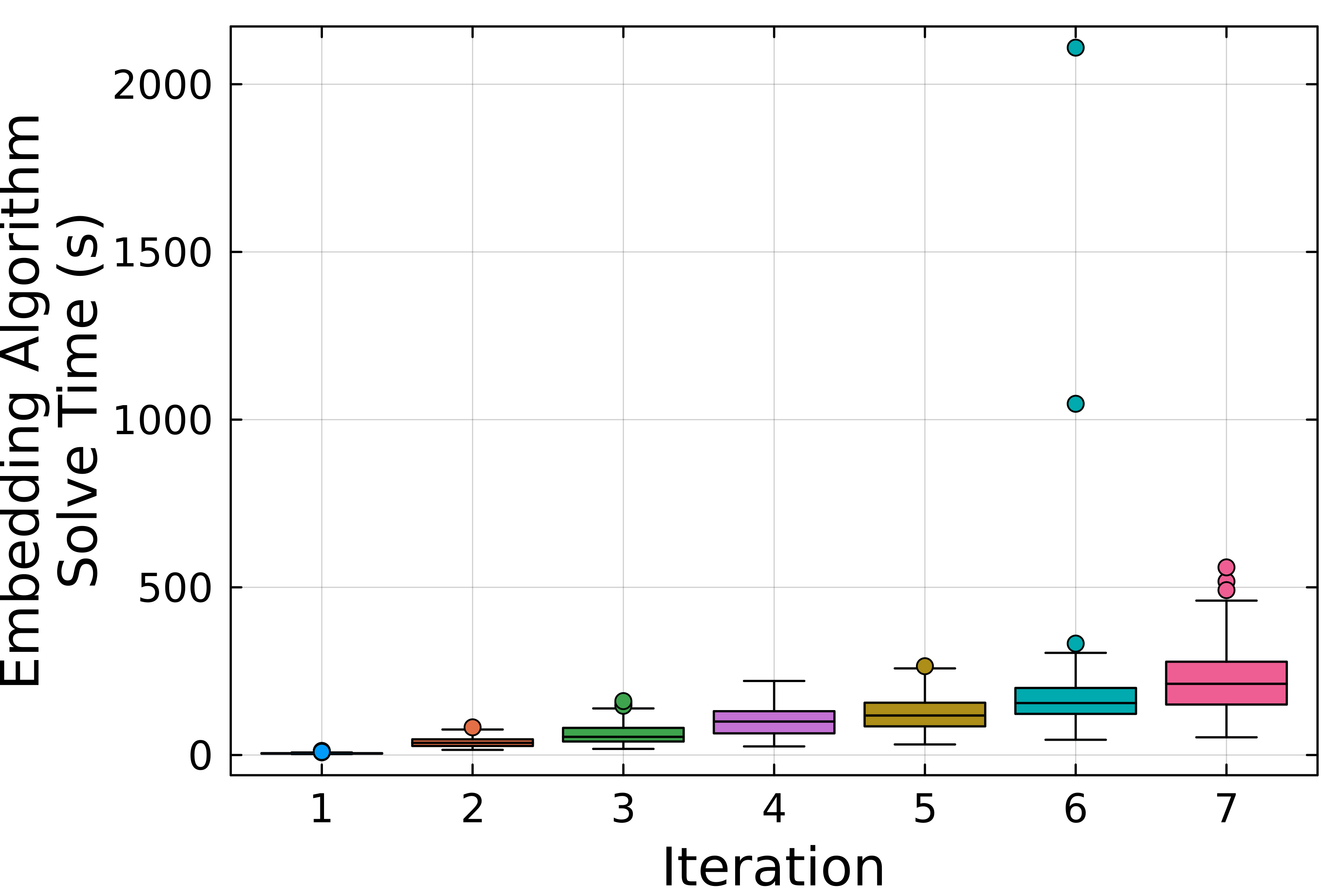}
    \caption{Embedding algorithm solve time for QC4UC with QA for 100 independent tests.}
    \label{sec:results:fig:embed-time}
\end{figure}
Figure \ref{sec:results:fig:embed-time} plots the embedding algorithm solve time, which is the total amount of time that it takes for the QA sampler to find an embedding for the problem using the methodology described in Cai et al.'s work \cite{APrctclHeurstcForFndngGrphMnrs2014}. 
The results show that the embedding algorithm solve time increases for each iteration of the algorithm. 
As the problem size of the QUBO MP increases with each cut, the computational complexity of the embedding search problem increases.
These results highlight the importance of qubit reduction techniques, such as the proposed rounded cuts and dynamic encoding regime implemented in QC4UC.

\begin{figure}[ht]
    \centering
    \subfigure[]{
  \includegraphics[width=0.47\linewidth]{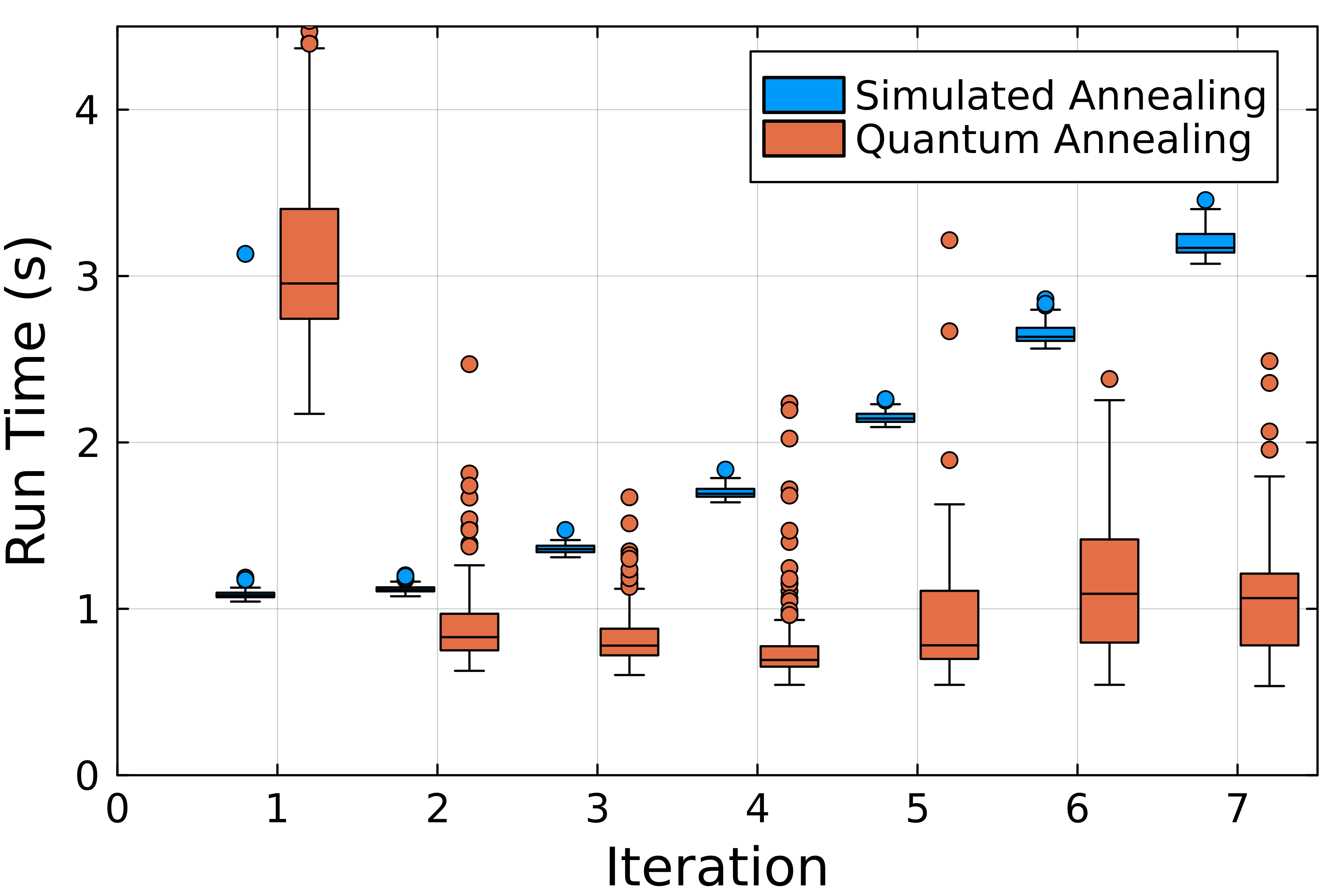}
  }
  \hspace{0.05cm}
    \subfigure[]{
  \includegraphics[width=0.47\linewidth]{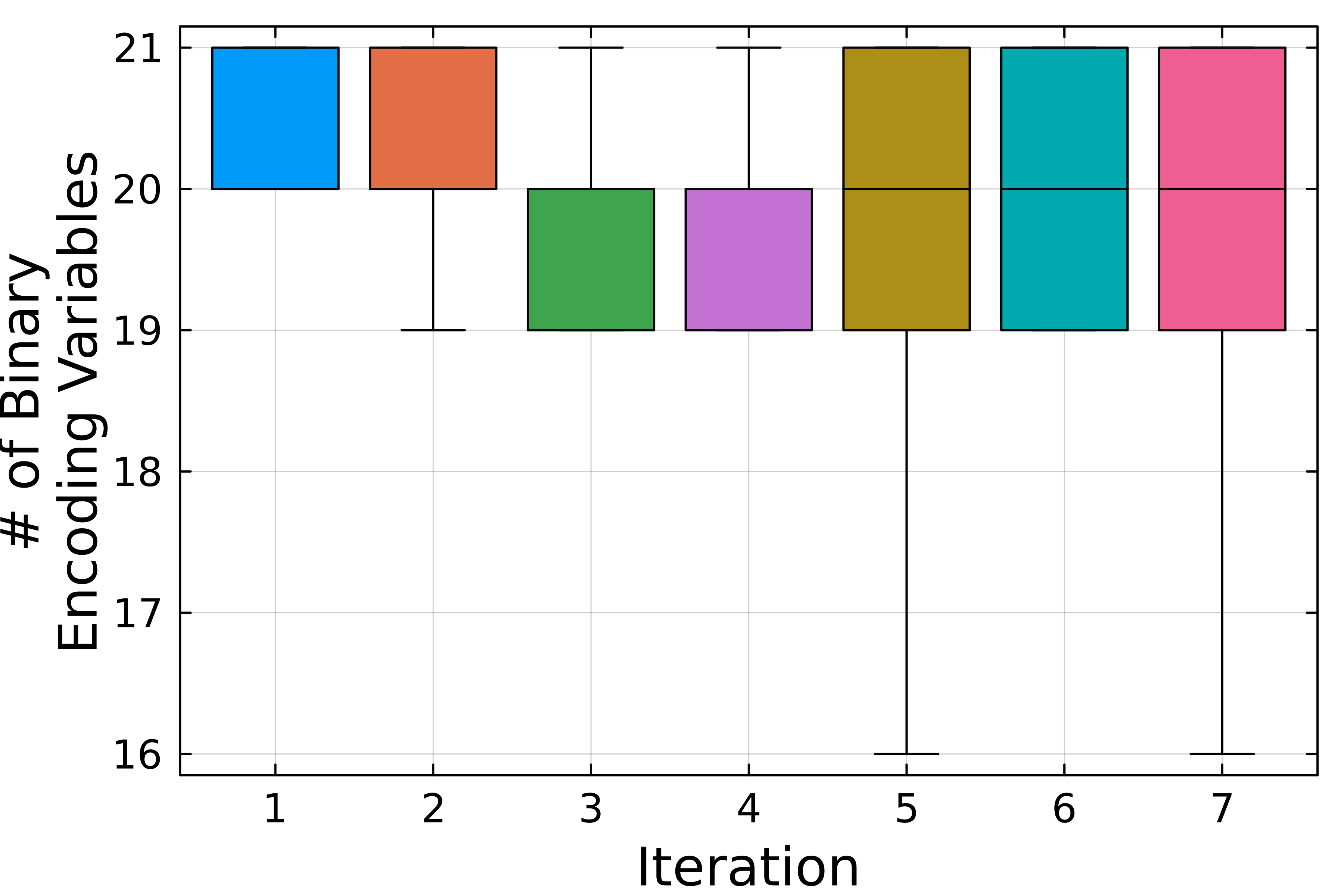}
  }
    \caption{(a) Total QC4UC run time per iteration, including MP solve time, SP solve time, and problem build time for QC4UC with SA and QC4UC with QA, excluding embedding algorithm solve time for QC4UC with QA. (b) Number of binary encoding variables used to approximate \(\boldsymbol{\eta}\) at each iteration of QC4UC with QA.}
    \label{sec:results:fig:solve-time_eta-encoding}
\end{figure}

Outside of the embedding algorithm solve time, the total algorithm run time for each iteration of QC4UC with QA demonstrates the linear scaling advantage of QAH, as shown in Figure \ref{sec:results:fig:solve-time_eta-encoding}.
The total algorithm run time includes the time to solve the MP, time to solve the SP, and, when using QA, the actual embedding time.
As the problem size grows with each iteration, QC4UC with QA scales better than QC4UC with SA.
It is seen clearly that, despite the stochasticity of the QA solver, QC4UC with QA can solve the larger problems faster on average than QC4UC with SA.
This demonstrates the necessity to test quantum algorithms on real quantum hardware, as classical computers cannot physically emulate the time scaling advantages of QA.
It is important to note that the first iteration of QC4UC with QA has a much longer run time than the rest of the iterations.
This is a result of the embedding of the QUBO into the QAH.
In the first iteration, the full problem must be embedded, whereas the subsequent iterations must only add the new cuts to the embedding.

\subsection{Dynamic Precision Encoding}
From Figure \ref{sec:results:fig:solve-time_eta-encoding}, we can see that the dynamic encoding of $\boldsymbol{\eta}$ allows us to reduce the qubit requirement over a static $\boldsymbol{\eta}$ encoding regime.
In these tests, QC4UC solved within 7 to 8 iterations and varied between 16 and 21 qubits per iteration for encoding the estimated secondary cost.
Without rounded cuts, 60 to 83 qubits were used per iteration to estimate the secondary cost.
This total qubit savings represents a significant portion of the available qubits for currently available QAH.
As QAH continues to increase in size, this dynamic encoding step will allow the solution of larger and more realistic systems sooner than traditional static encoding regimes.

These results showcase the ability of QC4UC to enhance the solution accuracy of QA while minimizing qubit requirements. 
This is in contrast to results in prior literature, which demonstrate algorithmic performance using SA rather than QA or used small test cases. 
Our results have demonstrated remarkable ability to achieve high-quality results for large-scale problems on current NISQ-era hardware.
QC4UC is able to reduce qubit requirements for cuts by about one third per iteration and achieve near-optimal solutions in about 50\% of cases.
Neither of these results have been seen before in the literature for a large-scale, complex UC instance like the one presented here.

\section{Conclusion}
The novel QC4UC algorithm introduces a rounded Benders' optimization cut formulation that significantly reduces binary precision encoding requirements compared to previous HQC algorithms.
This improvement allows QC4UC to significantly lower qubit requirements, enabling the solution of larger UC instances on real NISQ hardware than any previous work, to the best of our knowledge.
The QA-enabled algorithm achieves objective values within 1\% of the Gurobi optimal solution in 51\% of tests, showcasing its effectiveness in handling moderately sized power systems optimization problems.
A key feature of this algorithm is the novel \(k\)-local search recovery step, which addresses noisy or near-feasible solutions. 
It refines suboptimal or infeasible solutions caused by the inherent QA noise, enhancing overall solution quality.
However, a trade-off between solution quality and qubit requirements is observed, as expected in HQC approaches on NISQ hardware. 
This is addressed using a dynamic binary encoding step, enabling QC4UC to reduce qubit requirements where possible.
This dynamic binary encoding approach also applies to broader MIPs using Benders' decomposition, enabling research on larger problems with real NISQ hardware.

The current method significantly improves scalability; however, further research could focus on several areas. First, reducing qubit requirements while maintaining tighter bounds on solution quality by refining the rounding technique or exploring alternative encoding strategies. Second, embedding algorithm solve time can be substantially reduced by dynamic cut management techniques. Finally, incorporating more complex and realistic UC formulations, including non-convex power flow physics and combined-cycle generators.

\begin{credits}
\subsubsection{\ackname}
The authors acknowledge funding from the U.S. Department of Energy (DOE) under the Advanced Grid Modeling (AGM) program for the project "Rigorous Assessment and Development of Quantum and Hybrid Quantum-Classical Computing Approaches for Solving the Unit Commitment Problem." They also acknowledge the use of resources provided by the LANL's Institutional Computing Program, supported by the DOE National Nuclear Security Administration under Contract No. 89233218CNA000001. The authors would like to acknowledge Dr. Rabab Haider of University of Michigan Ann Arbor for her helpful feedback.

\subsubsection{\discintname}
The authors have no competing interests to declare that are relevant to the content of this article.
\end{credits}
%
%
%
\bibliographystyle{splncs04}
\bibliography{bib/AI4Opt-Rosemary}

%
\end{document}